\documentclass[letterpaper,12pt,reqno,twoside]{amsart}

\usepackage{amssymb}
\usepackage{amsmath}
\usepackage[dvips]{graphicx}
\usepackage{psfrag}

\newcommand{\be}{\begin{eqnarray}}
\newcommand{\ee}{\end{eqnarray}}

\newcommand{\beq}{\begin{equation}}
\newcommand{\eeq}{\end{equation}}
\newcommand{\beqn}{\begin{equation*}}
\newcommand{\eeqn}{\end{equation*}}

\newcommand{\cfg}[1]{C_{f,g}(#1)}
\newcommand{\fhat}{\hat{f}}
\newcommand{\ghat}{\hat{g}}

\newcommand{\torus}{{\mathbb{T}^2}}
\newcommand{\R}{\mathbb{R}}
\newcommand{\Z}{\mathbb{Z}}
\newcommand{\C}{\mathbb{C}}
\newcommand{\N}{\mathbb{N}}
\newcommand{\defas}{\mathrel{\raise.095ex\hbox{$:$}\mkern-4.2mu=}}
\newcommand{\defasr}{\mathrel{=\mkern-4.2mu\raise.095ex\hbox{$:$}}}
\newcommand{\nonzero}{\setminus\{0\}}
\newcommand{\abs}[1]{\lvert#1\rvert}
\newcommand{\norm}[1]{\lVert#1\rVert}
\newcommand{\bnorm}[1]{{\lVert#1\rVert}_\beta}
\newcommand{\Ccone}{\mathcal{C}}
\newcommand{\Econe}{\mathcal{E}}
\newcommand{\inv}{^{-1}}

\newcommand{\ie}{\textit{i.e.}}

\DeclareMathOperator{\tr}{tr}
\DeclareMathOperator{\Prob}{\mathbb{P}}

\newtheorem{theo}{Theorem}
\newtheorem{cor}[theo]{Corollary}
\newtheorem{lem}[theo]{Lemma}
\newtheorem{defn}{Definition}

\newtheorem{remark}[theo]{Remark}

\begin{document}

\title[Randomly Chosen Hyperbolic Toral Automorphisms]{Exponential Decay of Correlations for Randomly Chosen Hyperbolic Toral Automorphisms}
\author[Arvind Ayyer]{Arvind Ayyer$^1$}
\author[Mikko Stenlund]{Mikko Stenlund$^2$}
\address{${}^1$ Department of Physics \\
Rutgers University\\
136 Frelinghuysen Road\\
Piscataway, NJ 08854, USA}
\address{${}^2$ Department of Mathematics \\
    Rutgers University \\
    110 Frelinghuysen Road \\
    Piscataway, NJ 08854, USA}
\email{\{ayyer,mstenlun\}@math.rutgers.edu}


\date{\today}
\thanks{M. S. would like to thank the Finnish Cultural Foundation for funding. The work was supported in part by NSF DMR-01-279-26 and AFOSR AF 49620-01-1-0154.}

\begin{abstract}
We consider pairs of toral automorphisms $(A,B)$ satisfying an invariant cone property. At each iteration, $A$ acts with probability $p\in (0,1)$ and $B$ with probability $1-p$. We prove exponential decay of correlations for a class of H\"older continuous observables. 
\end{abstract}

\maketitle

\section{Introduction}

\subsection{Background}
Toral automorphisms are the simplest examples of Anosov maps. For deterministic Anosov maps, many ergodic and statistical properties such as ergodicity, existence of SRB measures and exponential decay of correlations for H\"older continuous observables are known.

Dynamical systems with randomness have been studied extensively in recent years. A typical model has been an Anosov map with noise. There are several good books on the subject \cite{Kifer,Arnold}. 

Products of random matrices are used in physics to model magnetic systems with random interactions and localization of electronic wave functions in random potentials.  They also play a central role in chaotic dynamical systems. In such applications, Lyapunov exponents provide information on the thermodynamic properties, electronic transport, and sensitivity for initial conditions. For more on the applications of products of random matrices, see \cite{Crisanti-Paladin-Vulpiani}.

We consider the action of two separate toral automorphisms $A$ and $B$, satisfying a cone condition, but which cannot be considered perturbations of one another. At each iteration, the matrix $A$ is picked with a certain probability and the matrix $B$ is picked if $A$ is not picked and applied on the torus. 

We are interested in the ergodic and statistical properties of the model for fixed realizations of the sequence of $A$'s and $B$'s obtained in this way. One might consider the ``environment'' (given by the sequence of matrices) to be fixed and the randomness to be associated with choosing the initial point on the torus. That is, our point of view is quenched randomness; one rolls the dice and lives with the outcome. 


\subsection*{Acknowledgements}
We are grateful to Giovanni Gallavotti, Sheldon Goldstein, Joel Lebowitz, Carlangelo Liverani, David Ruelle, and Lai-Sang Young for useful discussions.

\subsection{Toral automorphisms} 
Let $\torus$ be the 2-torus $\R^2/2\pi\Z^2$.

\begin{defn}
A map $A : \torus \circlearrowleft$  by the matrix action $x \mapsto Ax \pmod{2\pi}$ is called a \emph{toral automorphism} if the matrix $A$ has integer entries and $\det A = \pm 1$. It is a \emph{hyperbolic} toral automorphism if, further, the eigenvalues of the matrix $A$ have modulus different from 1.
\end{defn}

Since the eigenvalues of a $2\times 2$ matrix $A$ are given by the formula
\beq\label{eq:eigenvalues}
\frac{\tr A\pm \sqrt{(\tr A)^2-4\det A}}{2},
\eeq
we see that a toral automorphism is hyperbolic precisely when the eigenvalues are in $\R\setminus\{1\}$. The hyperbolicity condition reduces to 
\beq\label{eq:hyperbolic}
\begin{cases}
\abs{\tr A} > 2 & \text{if $\det A=+1$}, \\
\tr A\neq 0 & \text{if $\det A=-1$}.
\end{cases}
\eeq
Under the hyperbolicity assumption, the matrix has an eigenvalue whose absolute value is greater than 1, which we call the \emph{unstable eigenvalue} and denote by $\lambda^A_u$. Similarly, it has a \emph{stable eigenvalue}, $\lambda_s^A$, with absolute value less than 1. The corresponding eigenvectors $e_u^A$ and $e_s^A$ span linear subspaces $E_u^A$ and $E_s^A$, respectively. These we refer to as the \emph{unstable (eigen)direction} and the \emph{stable (eigen)direction}, respectively.

From now on, we will always assume that our toral automorphisms are hyperbolic and have determinant $+1$. This is necessary for the cone property formulated below.

It is important to notice that the eigenvalues in \eqref{eq:eigenvalues} are irrational. Consequently, the corresponding eigendirections $E_{s,u}^A$ have irrational slopes; if $A=\bigl(\begin{smallmatrix} a & b \\ c & d \end{smallmatrix}\bigr)$, the eigenvectors are given by the formula
\beq\label{eq:eigenvec}
e_{s,u}^A = \begin{pmatrix}
1 \\
(\lambda_{s,u}^A-a)/b
\end{pmatrix}.
\eeq
Here $b\neq 0$ by hyperbolicity; see \eqref{eq:hyperbolic}. Irrationality is equivalent to \beq\label{eq:irrational}
E_{u,s}^A \cap \Z^2 = \{0\}.
\eeq
An alternative, ``dynamical", way of proving \eqref{eq:irrational} starts with assuming the opposite; suppose that $0\neq q\in E_{s}^A \cap \Z^2$. Then $A^n q\to 0$ as $n\to\infty$, which is a contradiction because $A$ is an invertible integer matrix. The case of $E_u^A$ is similar using $A\inv$. 

In fact, the slopes of $E_{s,u}^A$---call them $\alpha_{s,u}^A$---are not only irrational but satisfy a stronger arithmetic property called the Diophantine condition: there exist $\epsilon>0$ and $K_\epsilon>0$ such that 
\beq
(q_1,q_2)\in \Z^2, \, q_1\neq 0\quad \Longrightarrow \quad \Bigl | \alpha_{s,u}^A-\frac{q_2}{q_1}  \Bigr | \geq \frac{K_\epsilon}{\abs{q_1}^{2+\epsilon}}.
\eeq
This tells us that in order for $E_{u,s}^A$ to come close to a point on the integer lattice $\Z^2\nonzero$, that point has to reside far away from the origin.

Finally, we point out that a toral automorphism is symplectic. That is, setting
\beq
J\defas \begin{pmatrix} 0 & 1 \\ -1 & 0 \end{pmatrix},
\eeq
any $2\times 2$ matrix $A$ with determinant one satisfies
\beq
A^T J A = J.
\eeq
For future reference, we define
\beq\label{eq:Atilde}
\widetilde A \defas (A^T)\inv = J A J\inv.
\eeq

\subsection{Invariant cones}
\begin{figure}[h!] 
\psfrag{EC}{$\Econe$}
\psfrag{CC}{$\Ccone$} 
\psfrag{eca}{$e^A_u$}
\psfrag{ecb}{$e^B_u$}
\psfrag{cca}{$e^A_s$}
\psfrag{ccb}{$e^B_s$}
\includegraphics[width=10cm,angle=0]{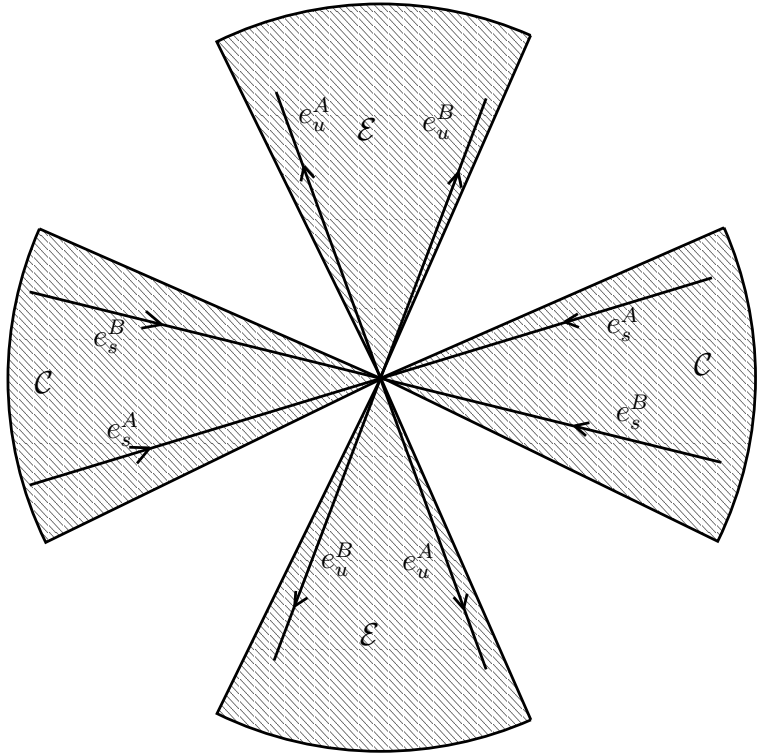}
\caption{An example of expansion and contraction cones.}
\label{fig:cones}
\end{figure}

A matrix on $\R^2$ maps lines running through the origin into lines running through the origin. Therefore, it is natural to consider \emph{cones}, \ie, sets in $\R^2$ spanned by two lines intersecting at the origin \cite{Alekseev,Wojtkowski}. Cones have also been used in the study of the spectrum of the transfer operator. See, for example, \cite{Blank-Keller-Liverani,Gouezel-Liverani,Baladi}.
\begin{defn}[Cone property]
A pair $(A_0,A_1)$ of hyperbolic toral automorphisms has the cone property if the following cones exist (Fig.~\ref{fig:cones}):
\\
An \emph{expansion cone}, $\Econe$, is a cone such that
\begin{enumerate}
\item $A_i\Econe\subset\Econe$,
\item there exists $\lambda_\Econe>1$ such that $\abs{A_i x}\geq \lambda_\Econe\abs{x}$ for $x\in\Econe$,
\item The $E_u^{A_i}$ do not lie along the boundary $\partial\Econe$: $E_u^{A_i}\cap \partial\Econe=\{0\}$.
\end{enumerate}
A \emph{contraction cone}, $\Ccone$, is a cone such that $\Econe\cap\Ccone=\{0\}$ and
\begin{enumerate}
\item $A_i\inv\Ccone\subset\Ccone$,
\item there exists $\lambda_\Ccone<1$ such that $\abs{A_i\inv x}\geq \lambda_\Ccone\inv \abs{x}$ for $x\in\Ccone$,
\item The $E_s^{A_i}$ do not lie along the boundary $\partial\Ccone$: $E_s^{A_i}\cap \partial\Ccone=\{0\}$.
\end{enumerate}
\label{defn:coneprop}
\end{defn}

\begin{remark}\label{rem:cones} 
A cone is a contraction cone with rate $\lambda_\Ccone$ (respectively expansion cone with rate $\lambda_\Econe$) for $(A,B)$ if and only if it is an expansion cone with rate $\lambda_\Ccone\inv$ (respectively contraction cone with rate $\lambda_\Econe\inv$) for $(A\inv,B\inv)$. With the aid of \eqref{eq:Atilde}, one checks that if one of the pairs $(A,B),(A^T,B^T),(A\inv,B\inv),(\widetilde A,\widetilde B)$ has the cone property then all of them do (with different cones). Moreover, the rates coincide for the corresponding cones of $(A,B)$ and $(\widetilde A,\widetilde B)$.
\end{remark}

The name ``expansion cone" is obvious, whereas ``contraction cone" deserves some caution: 
$
\abs{A_i x}\leq \lambda_\Ccone\abs{x}
$
holds in general \emph{only} under the assumption $A_i x\in\Ccone$, as opposed to the weaker $x\in\Ccone$.

Given a line $L$ passing through the origin transversely to $E_s^{A_i}$, it is a consequence of hyperbolicity that the image line $A_i^nL$ tends to $E_u^{A_i}$ as $n\to\infty$. Therefore, an expansion cone has to contain the unstable eigendirection, $E_u^{A_i}$. Similarly, considering the backward iterates, a contraction cone has to contain the stable eigendirection, $E_s^{A_i}$. In brief,
\beq\label{eq:cone-eigendir}
E_u^{A_i}\subset \Econe \quad\text{and}\quad E_s^{A_i}\subset \Ccone.
\eeq
The expansion and contraction rates in the cones are naturally bounded by the eigenvalues of the matrices:
\beq\label{eq:cone-eigenval}
1 < \lambda_\Econe<\abs{\lambda_u^{A_i}} \quad\text{and}\quad 1>\lambda_\Ccone>\abs{\lambda_s^{A_i}}.
\eeq


In particular, our results will hold for any two toral automorphisms with positive entries, as well as for any two toral automorphisms with negative entries. This is because the union of the first and third quadrant is automatically an expansion cone and the complement a contraction cone, as can be easily checked. By inverting the matrices, our results apply just as well to any two toral automorphisms whose diagonal elements are positive (respectively negative) and off-diagonal ones negative (respectively positive).  

\subsection{Random toral automorphisms}
For definiteness, let $A$ be chosen with probability $p$ and $B$ be chosen whenever $A$ is not chosen.

In order to model randomness, we first set
\beq 
A_0\defas A \quad\text{and}\quad A_1\defas B.
\eeq
On the space of sequences, $\Omega\defas\{0,1\}^\N$, we define the shift operator
\be
\tau:(\omega(0),\omega(1),\dots) \mapsto (\omega(1),\omega(2),\dots),
\ee
and independently for each index $n\in\N$ prescribe the probability $p$ to ``$\omega(n)=0$" and the probability $1-p$ to ``$\omega(n)=1$". The resulting product measure, $\Prob$, is a $\tau$-invariant ergodic probability measure on $\Omega$. 

The map
\beq\label{eq:skew}
\Phi:\Omega\times\torus \circlearrowleft\; : (\omega,x)\mapsto (\tau\omega,A_{\omega(0)}x)
\eeq
is called a skew product and defines a random dynamical system. If ${\mathrm m}$ stands for the normalized Lebesgue measure on $\torus$, \ie,
\beq
d{\mathrm m}=\frac{dx}{(2\pi)^2},
\eeq
then $\Prob\times {\mathrm m}$ is a $\Phi$-invariant probability measure, because ${\mathrm m}$ is $A_i$-invariant for $i=0,1$. As usual, we will write $\mu(f)\defas \int d\mu f$ for the integral of a function $f$ over a measure space with some measure $\mu$.

A basis of the space $L^2(\Omega,\Prob)$ of square integrable functions on $\Omega$ can be constructed as follows. Set $\sigma_i(\omega)\defas \sqrt{p/(1-p)}$ if $\omega(i)=1$, and $\sigma_i(\omega)\defas -\sqrt{(1-p)/p}$ if $\omega(i)=0$. Then define $\sigma_A \defas \prod_{i\in A}\sigma_i$ for any finite subset $A$ of $\N$. The set $\{\sigma_A\,|\, \text{$A\subset \N$ finite}\}$ is a countable orthonormal basis of $L^2(\Omega,\Prob)$.

Let us denote
$
A_\omega^n \defas A_{\omega(n-1)}\cdots A_{\omega(0)},
$
such that 
$
A_\omega^n = A_{\tau^k\omega}^{n-k}\, A_\omega^k 
$
for $1\leq k < n$.
Using the notation $\widetilde A \defas (A^T)\inv$ for any matrix $A$, we get
\beq\label{eq:tilde}
\widetilde A_\omega^n \defas \widetilde A_{\omega(n-1)}\cdots \widetilde A_{\omega(0)} = ((A_\omega^n)^T)\inv. 
\eeq
Moreover,
\beq\label{eq:cocycle}
\widetilde A_\omega^n=\widetilde A_{\tau^k\omega}^{n-k}\, \widetilde A_\omega^k \quad\text{for}\quad k=1,\dots,n-1.
\eeq

\subsection{Sensitive dependence on initial conditions}
As an ingredient of chaos, we discuss how the iterates of a point on the torus depend sensitively on the point chosen. By this we roughly mean that the distance, $\abs{A_\omega^n x-A_\omega^n y}$, between iterates of two close by points, $x$ and $y$, typically diverges at an exponential rate as $n$ grows.

To this end, we use the Multiplicative Ergodic Theorem (MET) that is originally due to Oseledets in the context of differential equations and smooth flows; see \cite{Oseledets} and also \cite{RuelleLyapunov}. 
\begin{theo}[$2\times 2$ MET]
Suppose that $\omega\mapsto A_{\omega(0)}$ is a measurable mapping from $\Omega$ to the space of real $2\times 2$ matrices and the mapping $\omega\mapsto \ln^+\norm{A_{\omega(0)}}$ is in $L^1(\Omega,\Prob)$. Here $\ln^+t\equiv \max(\ln t,0)$ and  $\norm{\,\cdot\,}$ is any matrix norm. Then there exists a set $\Gamma\subset\Omega$ with $\tau\Gamma\subset \Gamma$ and $\Prob(\Gamma)=1$ such that the following holds if $\omega\in\Gamma$:
\begin{enumerate}

\item The limit
\beq
\lim_{n \to \infty} \left((A^n_\omega)^T A^n_\omega \right)^{1/2n} = \Lambda_\omega
\eeq
\noindent
exists.

\item Let $\exp \chi_\omega^{(1)} < \cdots < \exp \chi_\omega^{(s)}$ be the eigenvalues of $\Lambda_\omega$ with $U_\omega^{(1)},\cdots,U_\omega^{(s)}$ the corresponding eigenspaces. Further, denote
\be
V_\omega^{(0)} & = & \{0 \} \\
V_\omega^{(r)} & = & U_\omega^{(1)} \oplus \cdots \oplus U_\omega^{(r)} \quad \mathrm{for} \; r=1,\cdots,s
\ee
\noindent
Then $s$ is either $1$ or $2$ and for $r=1,\dots,s$ we have
\beq
\lim_{n\to\infty}\frac{1}{n} \ln\abs{A_\omega^n x} = \chi_\omega^{(r)} \quad\text{for}\quad x\in V_\omega^{(r)}\setminus V_\omega^{(r-1)}.
\eeq
\noindent
The numbers $\chi_\omega^{(r)}$ are called Lyapunov exponents. Both $s$ and $\chi_\omega^{(r)}$ generically depend on $\omega$ but are $\tau$-invariant.
\end{enumerate}
\end{theo}

Notice that $\det \Lambda_\omega=1$, or $\sum_{r=1}^s \chi_\omega^{(s)}=0$, yields  two possibilities:
\begin{enumerate}
\item $s=1$ and $\chi_\omega^{(1)}=0$,
\item $s=2$ and $\chi_\omega^{(2)}=-\chi_\omega^{(1)} > 0$.
\end{enumerate}
The MET guarantees that the Lyapunov exponents are invariant under the flow $\omega\mapsto\tau\omega$. By ergodicity of the $\tau$-invariant measure $\Prob$, they are constant almost surely ($\Prob=1$). 

\begin{cor}\label{cor:tilde-Lyapunov} 
Replacing $A_\omega^n$ by $\widetilde A_\omega^n$ (see \eqref{eq:tilde}) does not change the Lyapunov exponents. In fact, setting $J\defas \bigl(\begin{smallmatrix} 0 & 1 \\ -1 & 0 \end{smallmatrix}\bigr)$,
\beq
\lim_{n \to \infty} \left((\widetilde A^n_\omega)^T \widetilde A^n_\omega \right)^{1/2n} = J \Lambda_\omega J\inv=\widetilde\Lambda_\omega=(\Lambda_\omega)\inv.
\eeq
\end{cor}
\begin{proof}
By \eqref{eq:Atilde}, $\widetilde A^n_\omega=J A^n_\omega J\inv$. Moreover, $J^T=J\inv$ and the symmetric matrix $(A^n_\omega)^T A^n_\omega$ is diagonalizable, such that $\bigl((\widetilde A^n_\omega)^T \widetilde A^n_\omega \bigr)^{1/2n}$ reads $J O_n \, \Lambda_n^{1/2n} \, O_n\inv J\inv=J \bigl((A^n_\omega)^T  A^n_\omega \bigr)^{1/2n} J\inv$ for some matrix $O_n$ and diagonal matrix $\Lambda_n$.
\end{proof}

With respect to the Lebesgue measure on $\R^2$, almost all (a.a.) points $x$ belong to the set $V_\omega^{(s)}\setminus V_\omega^{(s-1)}$ corresponding to the largest Lyapunov exponent $\chi_\omega^{(s)}$. $\Prob$-almost surely the latter equals a constant $\chi^{(s)}$. Therefore, for a.a.\@ $\omega$, for a.a.\@ x,
\beq
\lim_{n\to\infty}\frac{1}{n} \ln\abs{A_\omega^n x} = \chi^{(s)}.
\eeq
\textit{A priori}, there might not exist a \emph{positive} Lyapunov exponent ($s=1$). This would rule out sensitive dependence on initial conditions in the meaning of the notion described in the beginning of the subsection.

In particular, it does not follow from the classical works of Furstenberg \cite{Furstenberg}, Kesten \cite{Furstenberg-Kesten}, and Virtser \cite{Virtser} that the largest Lyapunov exponent is positive, because the Bernoulli measure used to choose a matrix at each step is concentrated at two points, $A$ and $B$, on $SL(2,\R)$. 

\begin{theo}[The largest Lyapunov exponent is positive]\label{thm:Lyapunov} 
Suppose that $(A,B)$ has the cone property. Then there are two distinct Lyapunov exponents, and $\chi^{(2)}=-\chi^{(1)}$. In fact, $0<\ln\lambda_\Econe\leq \chi^{(2)}\leq \max_{i} \ln \abs{\lambda_u^{A_i}}$ and $\min_i\ln\abs{\lambda_s^{A_i}}\leq \chi^{(1)}\leq \ln\lambda_\Ccone<0$. In particular, $\chi^{(2)}\geq \ln \max(\lambda_\Ccone\inv,\lambda_\Econe)$.
\end{theo}

\begin{proof}
The largest Lyapunov exponent is positive, because the expansion cone has nonzero measure. We conclude that $s=2$ in the MET.

Consider the intersection $E_{\omega}\defas \bigcap_{k\geq 0}(A^k_{\omega})\inv\Ccone$ of preimages of the contraction cone $\Ccone$. It is a line inside $\Ccone$ whose forward iterates forever remain in $\Ccone$. That is, if $x\in E_\omega$, then $A_\omega^n x\in \Ccone$ such that $\abs{A_\omega^n x}\leq \lambda_\Ccone^n \abs{x}$ for all $n\geq 0$. We must have $V_\omega^{(1)}=E_\omega$. The construction of $E_\omega$ is similar to that of random stable manifolds \cite{Young}.

The bounds are now obvious.
\end{proof}

\subsection{Observables and a H\"older continuity condition}
\begin{defn}
We say that a function (``observable") $f:\Omega\times\torus\to\C$ satisfies the strong H\"older condition with exponent $\beta\in[0,1]$, if 
\beq\label{eq:holder}
\bnorm{f}\defas \sup_{\omega\in\Omega}\sum_{q\in\Z^2}\abs{\hat f(\omega,q)}\abs{q}^\beta <
\infty.
\eeq
Here $\hat{f}$ is the Fourier transform of $f$.
\end{defn}

Observe that if $f$ satisfies \eqref{eq:holder}, then
\beq\label{eq:holder-decay}
\sup_\omega\abs{\hat f(\omega,q)} \leq \bnorm{f}\abs{q}^{-\beta}, \quad q\in\Z^2\nonzero.
\eeq
Because $0\leq\beta\leq 1$, $\abs{e^{it}-1}/\abs{t}^\beta$ is uniformly bounded in $t\in\R$, and we see that \eqref{eq:holder} implies H\"older continuity of $f(\omega,\cdot)$ with exponent $\beta$:
\be
\abs{f(\omega,x+y)-f(\omega,x)}\leq\sum_{q\in\Z^2\nonzero} \abs{\hat f(\omega,q)}\abs{q\cdot y}^\beta\frac{\abs{e^{iq\cdot y}-1}}{\abs{q\cdot y}^\beta}\leq C\abs{y}^\beta
\ee
for all $x,y\in\torus$. The opposite is not true; hence the adjective ``strong".

\subsection{Decay of correlations and mixing}
We define the $n$th (time) correlation function of two observables $f$ and $g$ as
\beq\label{eq:correlation}
\cfg{\omega,n} \defas \int_{\torus} d{\mathrm m}(x)  \, f(\omega,A_\omega^n x) g(\omega,x)-{\mathrm m}({f(\omega,\cdot)}) \, {\mathrm m}({g(\omega,\cdot)}).
\eeq
\noindent
We also need the related
\beq
C_{f,g}^\Phi(\omega,n) \defas \int_{\torus} d{\mathrm m}(x)  \, (f\circ \Phi^n \cdot g)(\omega,x) -{\mathrm m}({f(\tau^n\omega,\cdot)}) \, {\mathrm m}({g(\omega,\cdot)})
\eeq
and
\beq
C_{f,g}^\Phi(n) \defas \int_{\Omega \times \torus} d(\Prob\times{\mathrm m})  \, (f\circ \Phi^n \cdot g)-(\Prob\times{\mathrm m})(f) \, (\Prob\times{\mathrm m})(g),
\eeq
where $\Phi$ refers to the skew product \eqref{eq:skew}. Our result is the following:
\begin{theo}[Decay of Correlations]\label{thm:corr-decay} 
Let the pair $(A,B)$ satisfy the cone property (\ref{defn:coneprop}); see Remark~\ref{rem:cones}. There exist $c>0$ and $\rho>0$ such that, if $f$ and $g$ are two observables satisfying the strong H\"older condition with exponent $\beta\in (0,1]$, then for all $\omega\in\Omega$,
\beq\label{eq:corr-decay}
\sup_\omega\abs{\cfg{\omega,n}}\leq c\,\bnorm{f}\bnorm{g} \, e^{-\rho \beta n}
\eeq
and
\beq\label{eq:corr-decay-skew}
\sup_\omega\abs{C_{f,g}^\Phi(\omega,n)}\leq c\,\bnorm{f}\bnorm{g} \, e^{-\rho \beta n}
\eeq
hold for all $n\in\N$. In fact, we can take $\rho=\ln\min(\lambda_\Ccone\inv,\lambda_\Econe)$.

For all $\epsilon>0$, there exists a constant $C(\epsilon)$ such that, for almost all $\omega$, the upper bounds above can be replaced by
\beq\label{eq:corr-decay-Lyapunov}
c\,\bnorm{f}\bnorm{g} \,C(\epsilon)^{-\beta} \, e^{-(\chi^{(2)}-\epsilon)\beta n},
\eeq
where $\chi^{(2)}\geq \ln \max(\lambda_\Ccone\inv,\lambda_\Econe)$ is the a.e.\@ constant, positive, Lyapunov exponent.
\end{theo}

\begin{remark}
Without additional information concerning the convergence rates of $\frac{1}{n} \ln\abs{\widetilde A_\omega^n x}$ to the corresponding Lyapunov exponents, we have no control over $C(\epsilon)$ beyond the fact that it is an increasing function of $\epsilon$.
\end{remark}

The proof of Theorem~\ref{thm:corr-decay} follows in Section~\ref{sec:proof}. 

If $f,g$ are trigonometric polynomials (finite linear combination of exponentials of the form  $e^{iq\cdot x}$), they satisfy \eqref{eq:holder} trivially and thus, also \eqref{eq:corr-decay}. Trigonometric polynomials form a countable basis of $L^2(\torus,{\mathrm m})$. This implies that
$
\lim_{n\to\infty}\int_{\torus} d{\mathrm m}(x)  \, F(A_\omega^n x) G(x)={\mathrm m}(F) \, {\mathrm m}(G)
$
for any functions $F,G\in L^2(\torus,{\mathrm m})$. We say that every \emph{fixed} realization of the random sequence $(A_{\omega(0)}, A_{\omega(1)}, \dots)$ of maps is mixing on $\torus$.

Similarly, one should interpret \eqref{eq:corr-decay-skew} as a mixing result for the skew product, keeping $\omega$ fixed.

\begin{remark} 
It is true that the positivity of the Lyapunov exponent $\chi^{(2)}$ is enough for mixing, even if the cone condition is not satisfied. However, we need the cone condition to
\begin{enumerate}
\item check that $\chi^{(2)}$ actually is positive, and
\item obtain estimates on correlation decay, \ie, on the mixing rate.
\end{enumerate}
\end{remark}

\begin{cor}
If $(A, B)$ satisfies the cone property, then the skew product $\Phi$, or the random dynamical system, is mixing. If, moreover, $f$ and $g$ satisfy the strong H\"older condition and ${\mathrm m}({f(\omega,\cdot)}) \equiv {\mathrm m}({g(\omega,\cdot)}) \equiv 0$, then $\abs{C_{f,g}^\Phi(n)}\leq c\,\bnorm{f}\bnorm{g} \,\min\bigl(e^{-\rho\beta n}, C(\epsilon)^{-\beta}e^{-(\chi^{(2)}-\epsilon)\beta n}\bigr)$.
\end{cor}
\begin{proof}
For $f,g\in L^2(\Omega\times\torus,\Prob\times {\mathrm m})$, the difference $C_{f,g}^\Phi(n)-\Prob(C_{f,g}^\Phi(\cdot,n))$ has the expression
\beq
 \int_\Omega d\Prob(\omega)\, {\mathrm m}({f(\tau^n\omega,\cdot)}) \, {\mathrm m}({g(\omega,\cdot)}) - (\Prob\times{\mathrm m})(f) \, (\Prob\times{\mathrm m})(g),
\eeq
and tends to zero, because $\tau$ is mixing. 
Since $\{\sigma_A\,|\, \text{$A\subset \N$ finite}\}$ is a countable basis of $L^2(\Omega,\Prob)$, we get $L^2(\Omega\times\torus,\Prob\times {\mathrm m}) \cong L^2(\Omega,\Prob) \otimes L^2(\torus,{\mathrm m})$. The functions $\sigma_A(\omega)e^{iq\cdot x}$ form a countable basis of the latter and trivially satisfy \eqref{eq:holder} and an estimate corresponding to \eqref{eq:corr-decay-skew}. Hence, for any $f,g\in L^2(\Omega\times\torus,\Prob\times {\mathrm m})$,
$
\lim_{n\to\infty} C_{f,g}^\Phi(n)=0,
$
such that $\Phi$ is mixing. The second claim follows from \eqref{eq:corr-decay-skew} and \eqref{eq:corr-decay-Lyapunov} because $C_{f,g}^\Phi(n)=\Prob(C_{f,g}^\Phi(\cdot,n))$ when ${\mathrm m}({f(\omega,\cdot)}) \equiv {\mathrm m}({g(\omega,\cdot)}) \equiv 0$.
\end{proof}

The correlation function in \eqref{eq:correlation} has the Fourier representation
\beq\label{eq:correlation-Fourier}
\cfg{\omega,n}
=\sum_{q \in \mathbb{Z}^2\nonzero}\fhat(\omega,-\widetilde A_\omega^n q) \, \ghat(\omega,q).
\eeq
The proof of Theorem~\ref{thm:corr-decay} is based on the decay \eqref{eq:holder-decay} of $\hat f(q),\hat g(q)$ with increasing $\abs{q}$ and on controlling $\abs{\widetilde A_\omega^n q}$ with lower bounds that are increasing in $n$ but not too heavily decreasing in $\abs{q}$, such that summability persists. Here we gain by working with the sequence $(\widetilde A_\omega^n)_{n\in\N}$ instead of  $((A_\omega^n)^T)_{n\in\N}$, because the former has a Markov property---not shared by the latter---due to the order in which the matrix factors are multiplied. More precisely, if $q^n\defas \widetilde A_\omega^n q$, then $q^{n+1}$ is completely determined by $q^n$ and $\omega(n)$ as opposed to the entire history $(\omega(k))_{0\leq k\leq n}$.

\subsection{Comments}
Notice that Theorem~\ref{thm:corr-decay} is not needed for Theorem~\ref{thm:Lyapunov}; ergodicity of the shift $\tau$ alone is relevant for Theorem~\ref{thm:Lyapunov}.

Without affecting the proofs much, the cone property can be weakened by relaxing, for instance in the case of the expansion cone, the assumption that \emph{every} iteration results in expansion. Assuming instead the existence of a number $N$ such that $\abs{A_\omega^N x}\geq \lambda_\Econe \abs{x}$ for all $\omega$ if $x\in\Econe$ is equally sufficient.

In fact, it is true that the system is mixing even when the cones are merely invariant without any contraction/expansion assumption. More precisely, if $A_i$ are hyperbolic toral automorphisms and there exist non-overlapping cones $\Econe$ and $\Ccone$ with the properties $A_i\Econe\subset\Econe$ and $A_i\inv\Ccone\subset\Ccone$, then mixing occurs. This is so, because all possible products $\widetilde A_\omega^n$ turn out to be hyperbolic. However, we have no control over the mixing rate (the speed at which correlations decay) in this case.

Our results remain valid for any number of automorphisms, if all of the matrices have a mutual contraction cone and a mutual expansion cone with contraction and expansion rates bounded away from 1.  

One can pursue a line of analysis different from ours by emphasizing randomness in $\omega$. For instance, fixing $x\in\torus$, the sequence ${(A_\omega^n x)}_{n\in\N}$ of random variables on $\Omega$ is a Markov chain. Furthermore, with a more probabilistic approach, it should be possible to do without the cone property (excluding extreme cases such as $B=A\inv$ with $p=\tfrac{1}{2}$) and still obtain results for almost all $\omega$.

\section{Proof of Theorem~\ref{thm:corr-decay}} \label{sec:proof}
Before entering the actual proof, we explain our method for the sake of motivation. Starting with \eqref{eq:correlation-Fourier}, we have
\beq\label{eq:corr-bound-init}
\abs{\cfg{\omega,n}} \leq \sum_{q \in \mathbb{Z}^2\nonzero} \abs{\fhat(\omega,-\widetilde A_\omega^n q) \, \ghat(\omega,q)}.
\eeq
The idea is to split $\Z^2\nonzero$ into suitable pieces using cones. By the cone property and Remark~\ref{rem:cones}, we have the following at our disposal:

\vspace{3mm}

\noindent
{\bf Contraction cone, $\Ccone$.} With some positive $\lambda_\Ccone<1$,
\beq\label{eq:Cbound}
\abs{\widetilde A_i q} \leq \lambda_\Ccone\abs{q} \quad\text{if}\quad \widetilde A_i q\in\Ccone, \, i\in\{0,1\}.
\eeq
Moreover, $E_s^{\widetilde A_i}\subset \Ccone$. The complement $\Ccone^c$ is invariant in the sense that
\beq
q\in\Ccone^c \quad \Longrightarrow \quad \widetilde A_i q\in\Ccone^c.
\eeq

\vspace{1mm}

\noindent
{\bf Expansion cone, $\Econe$.} With some $\lambda_\Econe>1$,
\beq\label{eq:Ebound}
\abs{\widetilde A_i q} \geq \lambda_\Econe\abs{q} \quad\text{for}\quad q\in\Econe, \, i\in\{0,1\}.
\eeq
Moreover, $E_u^{\widetilde A_i}\subset\Econe$. The cone $\Econe$ itself is invariant in the sense that
\beq
q\in\Econe \quad \Longrightarrow \quad \widetilde A_i q\in\Econe.
\eeq

\vspace{1mm}

\noindent
{\bf The complement $(\Ccone\cup\Econe)^c$ and the number $M$.}
There exists a positive integer $M$ depending on only the choice of cones $\Ccone,\Econe$ such that, for any ``random" sequence $\omega$,
\beq
q\in (\Ccone\cup\Econe)^c\quad \Longrightarrow \quad \widetilde A_\omega^M q\in\Econe.
\eeq

Let us define
\beq
\lambda\defas \max{(\lambda_\Ccone,\lambda_\Econe\inv)} < 1.
\eeq
It follows from the bound \eqref{eq:Cbound} in the definition of $\Ccone$ that, if \linebreak $q,\cdots,\widetilde A_\omega^{n} q\in \Ccone$ and $q\neq 0$, then
\beq\label{eq:Ccone-temp}
1\leq\abs{\widetilde A_\omega^n q} \leq \lambda^{n}\abs{q},
\eeq
where the first inequality is due to the fact that $\widetilde A_\omega^n q$ is a nonzero integer vector. But the right-hand side of \eqref{eq:Ccone-temp} tends to zero as $n\to\infty$, which is a contradiction unless eventually $\widetilde A_\omega^n q\notin \Ccone$. We conclude that there exists a unique integer $N_\omega(q)$, called the \emph{contraction time}, satisfying
\beq
q\in\Ccone\nonzero \quad\Longrightarrow\quad \widetilde A_\omega^{N_\omega(q)} q \in \Ccone \quad\text{and}\quad \widetilde A_\omega^{N_\omega(q)+1} q \in \Ccone^c.
\eeq

Let us list some bounds used in proving Theorem~\ref{thm:corr-decay}. The proof of this lemma is given at the end of the section. 
\begin{lem}\label{lem:q-bounds}
The contraction time, $N_\omega(q)$, obeys the $\omega$-independent bound
\beq\label{eq:Ctime-bound}
N_\omega(q)\leq \frac{\ln\abs{q}}{\ln\lambda\inv}.
\eeq
The following complementary bounds hold:
\begin{align}\label{eq:Ccone-bound<}
1\leq & \abs{\widetilde A_\omega^n q} \leq \lambda^n\abs{q}  & \text{for} \quad & q\in\Ccone\nonzero, \, n\leq N_\omega(q), \\
\label{eq:Ccone-bound>}
&\abs{\widetilde A_{\omega}^{n} q} \geq C\lambda^{-n}\abs{q}\inv  & \text{for} \quad & q\in\Ccone\nonzero, \, n>N_\omega(q), \\
\label{eq:Ccompl-bound}
&\abs{\widetilde A_\omega^{n} q} \geq C\lambda^{-n}\abs{q}  & \text{for} \quad & q\in\Ccone^c, \, n\in\N.
\end{align}
\end{lem} 
In fact, \eqref{eq:Ctime-bound} is never used but we state it for the sake of completeness.

We now have a natural way of decomposing the sum in \eqref{eq:corr-bound-init} for \emph{each fixed value of $n$}, by means of the disjoint partition
\beq
\Z^2\nonzero=\{q\in\Ccone\,|\,n\leq N_\omega(q)\} \cup \{q\in\Ccone\,|\,n> N_\omega(q)\} \cup \Ccone^c\nonzero.
\eeq
Namely, we can rearrange the series (of nonnegative terms) as
\beq\label{eq:corr-decomp}
\sum_{q\in\Z^2\nonzero} = \sum_{\substack{q\in\Ccone\nonzero: \\ n\leq N_\omega(q)}} + \sum_{\substack{q\in\Ccone\nonzero: \\ n> N_\omega(q)}} + \sum_{q\in\Ccone^c}.
\eeq
In the first series on the right-hand side we use \eqref{eq:Ccone-bound<}, \ie, the fact that $\abs{q}$ is large for large $n$. In the second and the third series it is $\abs{\widetilde A_\omega^n q}$ that is large for large $n$, by \eqref{eq:Ccone-bound>} and \eqref{eq:Ccompl-bound}, respectively. 

Similarly, the last part of Theorem~\ref{thm:corr-decay} is based on the following refinement of Lemma~\ref{lem:q-bounds}:
\begin{lem}\label{lem:q-bounds-Lyapunov}
For all $\epsilon>0$, there exists a constant $C(\epsilon)$ such that, for almost all $\omega$, the bounds
\begin{align}\label{eq:XCcone-bound<}
1\leq & \abs{\widetilde A_\omega^n q} \leq C(\epsilon)\inv e^{-n(\chi^{(2)}-\epsilon)}\abs{q}  & \text{for} \quad & q\in\Ccone\nonzero, \, n\leq N_\omega(q), \\
\label{eq:XCcone-bound>}
&\abs{\widetilde A_{\omega}^{n} q} \geq C(\epsilon) e^{n(\chi^{(2)}-\epsilon)}\abs{q}\inv  & \text{for} \quad & q\in\Ccone\nonzero, \, n>N_\omega(q), \\
\label{eq:XCcompl-bound}
&\abs{\widetilde A_\omega^{n} q} \geq C(\epsilon) e^{n(\chi^{(2)}-\epsilon)}\abs{q}  & \text{for} \quad & q\in\Ccone^c, \, n\in\N,
\end{align}
hold true.
\end{lem}

\begin{proof}[Proof of Theorem~\ref{thm:corr-decay}]
We first prove \eqref{eq:corr-decay} and leave it to the reader to check that the upper bounds below apply just as well to the case of \eqref{eq:corr-decay-skew}. By the same token, we keep the first argument $\omega$ of the observables implicit and just write $\ghat(q)$ instead of $\ghat(\omega,q)$ \textit{etc}.  
Let us proceed case by case in the decomposition \eqref{eq:corr-decomp} of \eqref{eq:corr-bound-init}:

\vspace{3mm}
\noindent
{\bf Case $q\in\Ccone\nonzero, \, n\leq N_\omega(q)$.} By \eqref{eq:holder-decay} and \eqref{eq:Ccone-bound<},
\beq
\abs{\hat g(q)} \leq \bnorm{g}\abs{q}^{-\beta} \leq \bnorm{g} \lambda^{\beta n},
\eeq
such that
\beq
\sum_{\substack{q\in\Ccone\nonzero: \\ n\leq N_\omega(q)}} \abs{\fhat(-\widetilde A_\omega^n q) \, \ghat(q)} \leq 
\bnorm{g} \lambda^{\beta n}\sum_{\substack{q\in\Ccone\nonzero: \\ n\leq N_\omega(q)}} \abs{\fhat(-\widetilde A_\omega^n q)} \leq 
\bnorm{f} \bnorm{g} \lambda^{\beta n}.
\eeq

\vspace{3mm}
\noindent
{\bf Case $q\in\Ccone\nonzero, \, n> N_\omega(q)$.} By \eqref{eq:holder-decay} and \eqref{eq:Ccone-bound>},
\beq
\abs{\fhat(-\widetilde A_\omega^n q)} \leq \bnorm{f} \abs{\widetilde A_\omega^n q}^{-\beta} \leq
\bnorm{f} C^{-\beta} \lambda^{\beta n} \abs{q}^\beta,
\eeq
and we have
\beq
\begin{split}
\sum_{\substack{q\in\Ccone\nonzero: \\ n> N_\omega(q)}}\abs{\fhat(-\widetilde A_\omega^n q) \, \ghat(q)}  & \leq
\bnorm{f} C^{-\beta} \lambda^{\beta n}\sum_{\substack{q\in\Ccone\nonzero: \\ n> N_\omega(q)}}\abs{\ghat(q)}\abs{q}^\beta \\ 
& \leq
\bnorm{f}  \bnorm{g} C^{-\beta} \lambda^{\beta n}.
\end{split}
\eeq

\vspace{3mm}
\noindent
{\bf Case $q\in\Ccone^c$.} By \eqref{eq:holder-decay} and \eqref{eq:Ccompl-bound},
\beq
\abs{\fhat(-\widetilde A_\omega^n q)} \leq \bnorm{f}\abs{\widetilde A_\omega^n q}^{-\beta} \leq
\bnorm{f} C^{-\beta} \lambda^{\beta n} \abs{q}^{-\beta},
\eeq
which implies
\beq
\sum_{q\in\Ccone^c}\abs{\fhat(-\widetilde A_\omega^n q) \, \ghat(q)} \leq
\bnorm{f} C^{-\beta} \lambda^{\beta n} \sum_{q\in\Ccone^c}\abs{\ghat(q)}\abs{q}^{-\beta} \leq
\bnorm{f} \bnorm{g} C^{-\beta} \lambda^{\beta n}.
\eeq

\vspace{3mm}

Choosing
\beq
c\geq 1+2 C^{-\beta} \quad\text{and}\quad \rho\defas -\ln\lambda,
\eeq
the bound \eqref{eq:corr-decay} follows.

By the same argument, the fact that the upper bound \eqref{eq:corr-decay-Lyapunov} applies is an immediate consequence of Lemma~\ref{lem:q-bounds-Lyapunov}.
\end{proof}

\begin{proof}[Proof of Lemma~\ref{lem:q-bounds}]
Equation~\eqref{eq:Ccone-bound<} repeats \eqref{eq:Ccone-temp} and the related discussion, whereas \eqref{eq:Ctime-bound} is just another way of writing it in the case $n=N_{\omega}(q)$.

By the definition of $\Econe$, 
\beq\label{eq:Econe-bound}
\abs{\widetilde A_\omega^n q} \geq \lambda^{-n}\abs{q} \quad\text{when}\quad q\in\Econe, \, n\in\N.
\eeq
For each index $i$ and any $q\in\Z^2\nonzero$, we also have the \textit{a priori} bound
\beq\label{eq:apriori-bound}
\abs{\widetilde A_i q}\geq \mu_i \abs{q} \quad\text{with}\quad \mu_i\defas \norm{\widetilde A_i\inv}\inv=\abs{\lambda_s^{A_i}}.
\eeq
Therefore, if $q\in (\Ccone\cup\Econe)^c$ and $m\geq 1$, \eqref{eq:cocycle} implies
\be
\abs{\widetilde A_\omega^{M+m} q} \geq \lambda^{-m}\abs{\widetilde A_{\omega}^{M} q} \geq C\lambda^{-M-m}\abs{q},
\ee
where $C\defas\min_\omega(\lambda^{M}\mu_{\omega{(0)}}\dots\mu_{\omega{(M-1)}})$. Observing $C<1$ and sacrificing \eqref{eq:Econe-bound}, the latter bound extends to the whole of $\Ccone^c$, yielding \eqref{eq:Ccompl-bound}.

We also need a lower bound on $\abs{\widetilde A_\omega^n q}$---one that grows exponentially with $n$ but does not decrease too much with $\abs{q}$---assuming that $q\in\Ccone$ and $n>N_\omega(q)$. To this end, we notice $\widetilde A_{\omega}^{N_\omega(q)} q\in \Ccone^c$ and compute 
\be
\abs{\widetilde A_\omega^n q}=\abs{\widetilde A_{\tau^{N_\omega(q)}\omega}^{n-N_\omega(q)}\, \widetilde A_{\omega}^{N_\omega(q)} q} \geq C\lambda^{N_\omega(q)-n} \abs{\widetilde A_{\omega}^{N_\omega(q)} q}.
\ee
From \eqref{eq:Ccone-bound<} wee see $\abs{\widetilde A_{\omega}^{N_\omega(q)} q}\geq \lambda^{-N_\omega(q)}\abs{q}\inv$, so that \eqref{eq:Ccone-bound>} follows.
\end{proof}

\begin{proof}[Proof of Lemma~\ref{lem:q-bounds-Lyapunov}]

Modifying the constant, \eqref{eq:XCcone-bound>} follows from \eqref{eq:XCcone-bound<} and \eqref{eq:XCcompl-bound}, just like \eqref{eq:Ccone-bound>} follows from \eqref{eq:Ccone-bound<} and \eqref{eq:Ccompl-bound} in the proof of Lemma~\ref{lem:q-bounds}. 

$\widetilde A_\omega^n\Ccone^c$ becomes a thinner and thinner cone inside $\Econe$ as $n$ increases. If $x\in \partial \Ccone$ (the worst case), then $\frac{1}{n}\ln(\abs{\widetilde A_\omega^n x}/\abs{x})$ tends to $\chi^{(2)}$. Thus, for each $\epsilon>0$ and $n\in\N$, we have $\abs{\widetilde A_\omega^n x}/\abs{x} >C(\epsilon) e^{n(\chi^{(2)}-\epsilon)}$ for some choice of $C(\epsilon)$, and \eqref{eq:XCcompl-bound} follows.
 
Next, we prove \eqref{eq:XCcone-bound<}. Let $e_s(\omega)$ be a unit vector spanning the random stable line $\bigcap_{k\geq 0}(\widetilde A^k_{\omega})\inv\Ccone$ (see the proof of Theorem~\ref{thm:Lyapunov}). Then $\widetilde A_\omega^n e_s({\omega})\in\Ccone$ for all $n$.  Let us also pick an arbitrary unit vector $v$ in $\Econe$.


By an elementary geometric argument, there exists a constant $K$ such that if $x\in\Ccone$, $y\in\Econe$, \emph{and} $x+y\in\Ccone$, then $\abs{y}\leq K\abs{x}$ and $\abs{x}\leq K\abs{x+y}$ (the worst case: $x,x+y\in\partial \Ccone$ and $y\in\partial \Econe$). Recall that we are assuming $\widetilde A_\omega^n q\in \Ccone$ for each $n=0,\dots,N_\omega(q)$. If we split $q=q_1e_s(\omega)+q_2 v$, then $\widetilde A_\omega^n q=q_1\widetilde A_\omega^n e_s(\omega)+q_2 \widetilde A_\omega^n v$, where the first term belongs to $\Ccone$ and the second to $\Econe$. We gather that $\abs{q_2 \widetilde A_\omega^n v}\leq K\abs{q_1\widetilde A_\omega^n e_s(\omega)}$.

Because $\frac{1}{n}\ln(\abs{\widetilde A_\omega^n e_s(\omega)})$ tends to $-\chi^{(2)}$, for all $\epsilon>0$ there exists a $C_\epsilon$ such that, for all $n$, $\abs{\widetilde A_\omega^n e_s(\omega)}\leq C_\epsilon e^{-n(\chi^{(2)}-\epsilon)}$. In conclusion, $\abs{\widetilde A_\omega^n q}\leq C_\epsilon(1+K)\abs{q_1} e^{-n(\chi^{(2)}-\epsilon)}$, and
 \eqref{eq:XCcone-bound<} follows from $\abs{q_1}\leq K\abs{q}$.
\end{proof}


\end{document}